\documentclass[10pt]{amsart} 

\usepackage{amsfonts,amsmath,amsthm, graphics, amssymb}

\title[R.F. Shamoyan]{On Bergman projections and sharp decomposition theorems in tubular and related domains in $C^n$}

\author[]{R.F. Shamoyan}

\begin{document}
\begin{abstract} The theory of analytic function spaces in very general tubular domains over symmetric cones is a relatively new interesting research area. Tube domains are very general and very complicated domains. Recently several new results in this research area were provided in papers of B. Sehba and his coauthors concerning Bergman type operators in such type complicated unbounded domains. In this note we expand their results to certain spaces of analytic functions in products of tube domains. We define new integral operators of Bergman type and new analytic mixed norm spaces in such type domains and products of tube domains and provide new results on boundedness of certain Bergman type operators. Our results may have various nice applications in this research area.
Our results with very similar proofs may be valid in Siegel domains of second type, in bounded symmetric domains and bounded strongly pseudoconvex domains with the smooth boundary, various matrix domains. We will add at the end of this note a new sharp decomposition theorem for Bergman space. Previously such type sharp decomposition theorems in analytic function spaces were provided by author in other domains. Our sharp result in Bergman type function spaces enlarge that list of previously known such type assertions in analytic function spaces of several variables. We finnaly pose in addition various interesting new problems related to this research area and moreover indicate also some concrete schemes for solutions of these problems. We also provide in the second part of this note many interesting short comments and remarks.\\
Keywords: analytic functions, Bergman projections, tubular domains over symmetric cones, Bergman  spaces, product domains, mixed norm function spaces, decomposition theorems, harmonic function spaces, polydisk.
\end{abstract}

\maketitle 
\section{\textbf {Introduction}}
The goal of this paper is to provide new results on boundedness of Bergman type projections in tubular domains over symmetric cones and in products of such type domains. Note that for particular values of parameters our results are well-known. Our results on Bergman type projections may have various interesting applications in complex function theory of several variables namely in spaces of analytic functions in tubular domains over symmetric cones. We provide first basic notations and definitions of complex function theory in tubular domains over symmetric cones which are needed for this paper. The theory of analytic spaces in tubular domains over symmetric cones is an active research area (see, for example, \cite{10}-\cite{16} and references there).

Let $T_\Omega = V + i\Omega$ be the tube domain over an irreducible symmetric cone $\Omega$ in the complexification $V^{\mathbb{C}}$ of an $n-$dimensional Euclidean space $V.$ $\mathcal{H}(T_{\Omega})$ denotes the space
of all holomorphic functions on $T_{\Omega}$. Following the notation of \cite{14} and \cite{7} we denote
the rank of the cone $\Omega$ by $r$ and by $\Delta$ the determinant function on $V$.

Letting $V = \mathbb{R}^n$, we have as an example of a symmetric cone on $\mathbb{R}^n$ the Lorentz cone
$\Lambda_n$ which is a rank 2 cone defined for $n \geq 3$ by
$$\Lambda_n = \{y\in\mathbb{R}^n: y_1^2-\ldots-y_n^2>0, y_1 > 0\}.$$

The determinant function in this case is given by the Lorentz form
$$\Delta(y) = y_1^2-\ldots-y_n^2$$
(see, for example, \cite{14}).

Also, if $t, k \in \mathbb{R}^r$, then $t < k$ means $t_j < k_j$ for all $1\le j\le r$.

For $\tau\in\mathbb{R}_+$ and the associated determinant function $\Delta(x)$ \cite{14} we set
$$
A_{\tau}^{\infty}(T_{\Omega}) = \left\{F\in\mathcal{H}(T_{\Omega}): ||F||_{A_{\tau}^{\infty}} = \sup\limits_{x+iy\in T_{\Omega}}|F(x+iy)|\Delta^{\tau}(y)<\infty\right\}.
$$

It can be checked that this is a Banach space.
For $1\le p,q < +\infty$ and $\nu\in\mathbb{R}$, and $\nu>-1$ we denote by $A_{\nu}^{p,q}(T_{\Omega})$ the mixed-norm
weighted Bergman space consisting of analytic functions $f$ in $T_{\Omega}$ such that
$$
||F||_{A_{\nu}^{p,q}} = \left(\int_{\Omega} \left(\int_V |F(x+iy)|^p\,dx\right)^{\frac{q}{p}}\Delta^{\nu}(y)\,dy \right)^{\frac{1}{q}}<\infty.
$$

It will be interesting to study also these analytic function spaces when $p$ or $q$ are equal to $\infty.$ Such spaces can be easily defined, as in case of simpler domains.

We can extend easily these $A_{\nu}^{p,q}$ analytic spaces to product domains in at least two ways. The study of such type analytic spaces is a new interesting problem.

This is a Banach space.
Replacing above $A$ by $L$ we will get as usual the corresponding larger space of
all measurable functions in tube over symmetric cone with the same quazinorm (see \cite{10}, \cite{7}).
It is known that the $A_{\nu}^{p,q}(T_{\Omega})$ space is nontrivial if and only if $\nu>-1$ (see \cite{14}, \cite{9}).
When $p = q$ we write (see \cite{14})
$$A_{\nu}^{p,q}(T_{\Omega}) = A_{\nu}^{p}(T_{\Omega}).$$

This is the classical weighted Bergman space with usual modification when $p = \infty$.

The (weighted) Bergman projection $P_{\nu}$ is the orthogonal projection from the Hilbert space $L_{\nu}^2(T_{\Omega})$
onto its closed subspace $A_{\nu}^2(T_{\Omega})$ and it is given by the following integral formula (see \cite{14})
$$
P_{\nu}f(z) = C_{\nu}\int_{T_{\Omega}}B_v(z,w)f(w)dV_{\nu}w,
$$
where $$B_v(z,w) = C_{\nu}\Delta^{\nu+\frac{n}{r}}((z-\overline{w})/i)$$
is the Bergman reproducing kernel for $A_{\nu}^2(T_{\Omega})$ (see \cite{14}, \cite{7}).

Here we used the notation $dV_{\nu}(w) = \Delta^{\nu-\frac{n}{r}}(v)dudv$. Below and here we use constantly
the following notations $w = u + iv \in T_{\Omega}$ and also $z = x + iy \in T_{\Omega}$.
Hence for any analytic function from $A_{\nu}^2(T_{\Omega})$ the following integral formula is valid (see also \cite{14})
$$	f(z)=C_{\nu}\int\limits_{T_{\Omega}}B_v(z,w)f(w)dV_{\nu}w.
$$

In this case sometimes below we say simply that the $f$ function allows Bergman representation via Bergman kernel with $\nu$ index. Note that these assertions have direct copies in simpler cases of analytic function spaces in unit disk, polydisk, unit ball, upperhalfspace $\mathbb{C}_+$ and in spaces of harmonic functions in the unit ball or upperhalfspace of Euclidean space $\mathbb{R}^n$. These classical facts are well-known and can be found, for example, in  \cite{11}, \cite{6} and in some items from references there.
Above and throughout the paper we write $\mathbb{C}$ (sometimes with indexes) to denote positive constants which might be different each time we see them (and even in a chain of inequalities), but are independent of the functions or variables being discussed.

In the second part of this paper we formulate new sharp decomposition theorem for Bergman $A^1_{\alpha}$ function spaces in tube domain and pose various new interesting problems for readers concerning same type results in harmonic function spaces in several variables in the unit ball and upper half space of $\mathbb{R}^n$, indicating shortly some ways for readers on how to solve such type problems. Previously such sharp results were provided by the author in Bergman type function spaces in the unit ball and in bounded strongly pseudoconvex domains with smooth boundary in his recent papers with E. Tomashevskaya (see, for example, \cite{20}-\cite{22}).
Such type sharp results may have various applications in function theory of several variables.

Rather transparent arguments which we found in proof of theorem 3 in tube domains over symmetric cones and earlier in bounded strongly  pseudoconvex domains shows that similar type sharp decomposition theorems for Bergman spaces may be valid also with very similar proofs in Siegel domains of second type, bounded symmetric domains and minimal bounded homogeneous domains.

Proofs of  all our sharp decomposition assertions are based on special integral representations, uniform estimates for Bergman spaces and Forelly-Rudin estimates which are available in many domains.

\section{\textbf {Basic preliminaries on symmetric cones \\and determinant function}}

We first shortly remind the readers some basic facts on symmetric cones (see  \cite{4}, \cite{1}). A subset $\Omega$ of $\mathbb{R}^{n}$ or $V$, so that $\operatorname{dim} V=n$ to be a cone if $\lambda x \in \Omega$, for all $x \in \Omega, \lambda>0$, if $\lambda x+\mu y \in \Omega$ for all $x, y \in \Omega, \lambda, \mu>0$ then it is convex. Let in addition $\Omega^{*}=\left\{y \in \mathbb{R}^{n}:(y / x)>0\right.$, for all $\left.x \in \bar{\Omega} \backslash\{0\}\right\}$ and $\Omega^{*}=\Omega$. This type open cone is selfdual ($\Omega^{*}$ is dual cone).

Let $G(\Omega)=\left\{g \in G l\left(\mathbb{R}^{n}\right): g \Omega=\Omega\right\},$ where $G l\left(\mathbb{R}^{n}\right)$ denotes the group of all linear invertible transformation of $\mathbb{R}^{n}$. If for all $x, y \in \Omega, y=g x$, for some $g \in G(\Omega)$ then our open convex cone $\Omega$ is homogeneous, if also $\Omega^{*}=\Omega$ then it is symmetric cone. These are one of the main objects of this paper.

If the equation $\Omega=\Omega_{1}+\Omega_{2}$ is not possible for each $V_{1} \subset \mathbb{R}^{n}, $ $V_{2} \subset \mathbb{R}^{n}$, then our cone is irreducible, here $V_{i} \neq \varnothing, $ $ i=1,2,$ ( $\Omega_{1}, \Omega_{2}$ are symmetric cones), where also $\Omega_{i} \subset V_{i}, i=1,2.$

We remind shortly the reader now basic facts on determinant $\Delta^{t}(\operatorname{Im} z), $ $ z \in \mathbb{C}^{n}$, $t \in(0, \infty)$. We fix $V$ -- a simple Euclidean Jordan algebra with rank $r$.

(a) A Jordan algebra $V$ over $R$ is said to be Euclidean if there exists a positive definite bilinear symmetric form on $V$ which is associative $(L(x) u / v)=(u, L(x) v)$, where ( $u, v$ ) is an inner product on $R^{n}$, for all $x, u, v \in V$;

(b) A Jordan algebra is simple if all it's ideals are trivial;

(c) We define rank of $V$.

If $x \in V, $ $m(x)=\min \{k>0:\left(l, x, x^{2}, \ldots, x^{k}\right)$ are linearly dependent$\}$, then $1 \leq m(x) \leq \operatorname{dim} V$ and $r=\max \{m(x): x \in V\}$, we say rank of $V$ is $r$.

According to spectral theorem if $V$ has rank $r$, then $x=\sum_{i=1}^{r} \lambda_{i} c_{i}, $ $ \lambda_{i} \in \mathbb{R}, $ $ c_{i}$ are elements of so called Jourdan frame, and $\left\{\lambda_{i}\right\}$ are determined uniquely by $x$ (with their multiplicities). We fix now a Peirce decomposition of $V=\oplus_{1 \leq i \leq j \leq r} V_{i j}$; we formally look at V as a space of symmetric matrices $\left(V_{i j}\right), V_{i i}=R c_{i}$, where $R$ is a special mapping (see \cite{7}), $V_{i j}=V\left(c_{i}, 1 / 2\right) \cap V\left(c_{i}, 1 / 2\right)=\left\{x \in V: c_{i} x=c_{j} x=\frac{x}{2}\right\}$, $i<j, \operatorname{dim} V_{i j}=d=2 \frac{n / r-1}{r-1}).$ We denote by $P_{i j}$ the orthogonal projection of $V$ onto $V_{i j}$ for $i \leq j$. Finally we denote by $\Delta_{j}(x), $ $j=1, \ldots, r$, the principal minors of $x \in V$ with respect to the fixed Jordan frame $\left\{c_{1}, \ldots, c_{r}\right\}$. That is $\Delta_{k}(x)$ is the determinant of the projection $P_{k} x$ of $x$ in the Jordan subalgebra $V^{(k)}=\oplus_{1 \leq i \leq j \leq r} V_{i j}$.

It is well-known that $\Omega=\left\{x \in V: \Delta_{k}(x)>0, k=1, \ldots, r\right\}$. We have also $\Delta_{k}(m x)=\Delta_{k}(x), $ $x \in V, $ $m \in \mathbb{Z}_{+}, m>0$. See other properties of $\Delta_{k}$ in \cite{1}.

We define $\Delta_{s}(x)=\prod_{j=1}^{r} \Delta_{j}^{s_{j}-s_{j}+1}(x)=\Delta_{1}^{s_{1}-s_{2}}(x) \ldots \Delta_{r}^{S_{r}}(x), $ $x \in \Omega, $ $s \in C^{r}$. We have that $\left|\Delta_{s}\right|=\Delta(\operatorname{Im} z)$ and $\Delta_{s} \sum_{i=1}^{r} a_{i} c_{i}=\prod_{i=1}^{r} a_{i}^{s_{i}}, $ $ a_{i}>0, $ $i=1, \ldots, r$.

To formulate our theorems we will need basic facts of theory of analytic function spaces in tubular domains over symmetric cones taken from  \cite{4}, \cite{1}. Let $d v(w)$ and $d v_{\alpha}(w)=\left[\Delta^{\alpha-\frac{n}{r}}(v)\right] d u d v, $ $\alpha>n / r-1$, be a standard Lebesque measure in tubular domains over symmetric cone $T_{\Omega}$ and weighted Lebesque measure in tube, $w=u+i v$. The weighted Bergman kernel $B_{\nu}$ of $T_{\Omega}$ is given as usual by
$$
B_{\nu}(w, z)=\left(d_{\nu}\right) \Delta\left(\frac{w-\bar{z}}{i}\right)^{-\nu-\frac{n}{r}},\; w, z \in T_{\Omega}, \nu \in \mathbb{R},
$$
is a Bergman constant, where
$$
d_{\nu}=\left(c_{\nu}^{-1}\right) \Gamma\left(\nu+\frac{n}{r}\right).
$$

Let $\Omega$ be an irreducible symmetric cone in the Euclidean space $V$, and $T_{\Omega}=V+i \Omega$ the corresponding tube domain in the complexified space $V^{\mathbb{C}}$. We shall note $n$ the dimension of $V$ and $r$ the rank of $\Omega$. Moreover, we shall denote by $(x \mid y)$ the scalar product in $V$, and by $\Delta$ the determinant function. For the description of such cones $\Omega$ in terms of Jordan algebras, one may use the book of Faraut and Koranyi \cite{7}. One may also have in mind the typical example that one obtains when $V$ is the space of real symmetric $r \times r$ matrices, and $\Omega$ is the cone of positive definitive matrices. In this example, the scalar product on $V$ is induced by the Hilbert-Schmidt norm of the matrices, and the determinant function is given by the determinant of the matrices. The rank is $r$, while the dimension is $\frac{r(r+1)}{2}$.

We shall also make use the $g$ eneralized wave operator on $V$, given by $\square=\Delta\left(\frac{1}{i} \frac{\partial}{\partial x}\right)$. This is a differential operator of degree $r$, defined by the equality
$$
\Delta\left(\frac{1}{i} \frac{\partial}{\partial x}\right)\left[e^{i(z \mid \zeta)}\right]=\Delta(\zeta) e^{i(x \mid \zeta)}, \zeta \in V.
$$

It's name is due to another fundamental example, given by the forward light cone in $\mathbb{R}^{n}$,
$$
\left\{x \in \mathbb{R}^{n} ; x_{1}>\sqrt{x_{2}^{2}+\cdots+x_{n}^{2}}\right\},
$$
which is of rank 2 . In this case, the determinant function is equal to
$$
\Delta(x)=x_{1}^{2}-x_{2}^{2}-\cdots-x_{n}^{2}.
$$

\section{\textbf{Preliminaries on geometry of tubular domains \\over symmetric cones, basic lemmas, new analytic spaces}}

In this section we will collect several very useful assertions from \cite{10}, \cite{7}--\cite{9} mainly
concerning so-called r-lattices that will be used rather often in all proofs of our sharp
embedding theorems below.

Let $T_{\Omega} \subset \mathbb{C}^n$ be a bounded tubular domains over symmetric cones in $\mathbb{C}^n$. We shall
use the following notations:
\begin{itemize}
	\item $\delta: T_{\Omega}\rightarrow\mathbb{R}^+$ will denote the determinant function from the boundary, that is $\delta(z) = \Delta(Im z)$. Let $d\nu_t(Z) = (\delta(z))^{t}d\nu(z), t>-1$;
	\item $\nu$ will be the Lebesgue measure on $T_{\Omega}$;
	\item $H(T_{\Omega})$ will denote the space of holomorphic function on $T_{\Omega}$, endowed with the topology of uniform convergence on compact subsets;
	\item $B: T_{\Omega}\times T_{\Omega}\rightarrow\mathbb{C}$ will be the Bergman kernel of $T_{\Omega}$. Note that if $B$ is kernel of type $t, t\in\mathbb{N}$, then $B^s$ is kernel of type $st, s\in\mathbb{N}, t\in\mathbb{N}$. This follows directly from definition (see \cite{10}, \cite{7}--\cite{9}, \cite{4}, \cite{5}). Note $B = B_{2n/r}$ (see \cite{10}, \cite{7}--\cite{9}, \cite{4}, \cite{5});
	\item given $r\in(0,\infty)$ and $z_0\in T_{\Omega}$, we shall denote by $B_{T_{\Omega}}(z_0, r)$ the Bergman ball.
\end{itemize}

See, for example, \cite{10}, \cite{7}--\cite{9}, \cite{4}, \cite{5}, for definitions, basic properties and applications to geometric function theory of the Bergman distance and \cite{10}, \cite{7}--\cite{9}, \cite{4}, \cite{5} for definitions and basic properties of the Bergman kernel. Let us now recall a number of vital results proved in $T_{\Omega}$. The first two give information about the shape of Bergman balls:

\textbf{Lemma 1} (see \cite{10}, \cite{7}--\cite{9}, \cite{4}, \cite{5}).  {\it Let $T_{\Omega}\subset\mathbb{C}^n$ be a bounded tubular domains over symmetric cones, and $r\in(0, \infty)$. Then $$\nu(B_{T_{\Omega}}(\cdot, t))\approx\delta^{2r/n}.$$}

\textbf{Lemma 2} (see \cite{10}, \cite{7}--\cite{9}, \cite{4}, \cite{5}).  {\it Let $T_{\Omega}\subset\mathbb{C}^n$ be a bounded tubular domains over symmetric cones. Then there is $C > 0$ such that $$\frac{C}{1-r}\delta(z_0)\le\delta(z)\le\frac{1-r}{C}\delta(z_0)$$ for all $r\in(0, \infty), z_0\in T_{\Omega}$ and $z\in B_{T_{\Omega}}(z_0, r)$.}

\textbf{Definition 1}. Let $T_{\Omega}\subset\mathbb{C}^n$ be a tubular domains over symmetric cones, and $r > 0$. An $r$-lattice in $T_{\Omega}$ is a sequence ${a_k}\subset T_{\Omega}$ such that $T_{\Omega} = \bigcup\limits_{k}B_{T_{\Omega}} (a_k, r)$ and there exists $m > 0$ such that any point in $T_{\Omega}$ belongs to at most $m$ balls of the form $B_{T_{\Omega}}(a_k, R)$, where $R = \frac12(1 + r)$. Note by Lemma 2, $$\nu_{\alpha}(B_{T_{\Omega}}(a_k, R)) = \int\limits_{B_{T_{\Omega}}(a_k, R)} \delta^{\alpha}(z)d\nu(z) = (\delta^{\alpha}(a_k))\nu(B_{T_{\Omega}}(a_k, R)), \alpha>-1.$$
The existence of $r$-lattice in tubular domains over symmetric cones is ensured by the following

\textbf{Lemma 3} (see \cite{10}, \cite{12}-\cite{15}, \cite{7}--\cite{9}, \cite{4}, \cite{5}).  {\it Let $T_{\Omega}\subset\mathbb{C}^n$ be a bounded tubular domains over symmetric cones. Then for every $r\in(0, \infty)$ there exists an $r$-lattice in $T_{\Omega}$, that is there exists $m\in\mathbb{N}$ and a sequence ${a_k}\subset T_{\Omega}$ of points such that $T_{\Omega} = \bigcup\limits_{k=0}^{\infty}B_{T_{\Omega}}(a_k, r)$ and no point of $T_{\Omega}$ belongs to more than $m$ of the balls $B_{T_{\Omega}}(a_k, R)$, where $R =\frac{1}{2}(1 + r)$.}

We will call $r$-lattice sometimes the family $B_{T_{\Omega}}(a_k; r)$. Dealing with $B$ Bergman kernel we always assume $|B(z;a_k)|\asymp|B(a_k;a_k)|$ for any $z\in B_{T_{\Omega}}(a_k; r)$, $r\in(0;\infty)$ (see  \cite{10}, \cite{7}--\cite{9}, \cite{4}, \cite{5}). Let $m=(2n/r)l, l\in\mathbb{N}$. Then $|B_m(z;a_k)|\asymp|B_m(a_k;a_k)|, z\in B_{T_{\Omega}}(a_k; r), r\in(0;\infty)$. This fact is crucial for embedding theorems in tubular domains over symmetric cones (see also \cite{3}).

\textbf{Lemma 4} (see \cite{12}--\cite{15}, \cite{7}, \cite{4}, \cite{5}). {\it Let $T_{\Omega}\subset\mathbb{C}^n$ be a tubular
domains over symmetric cones. Given $r\in(0;\infty)$, set $R = \frac12(1+r)\in(0;\infty)$. Then there exists a $C_r > 0$ depending on $r$ such that $$\forall z_0\in T_{\Omega},\ \forall z\in B_{T_{\Omega}}(z_0,r),\ \chi(z)\le\frac{C_r}{\nu(B_{T_{\Omega}}(z_0,r))}\int\limits_{B_{T_{\Omega}}}\chi\,d\nu$$ for every nonnegative plurisubharmonic function $\chi:T_{\Omega}\rightarrow\mathbb{R}^+.$}

\textbf{Lemma 5} (see \cite{12}--\cite{15}, \cite{7}, \cite{4}, \cite{5}).{\it
\begin{enumerate}
	\item Let $\lambda>\frac{n}{r}-1$ be fixed. Then $\Delta(y+y')\ge \Delta(y) \forall y,$ $y'\in\Omega$, \\$|\Delta^{-\lambda}(\frac{x+iy}{i})|\ge\Delta(y)^{-\lambda}$, $\forall x\in\mathbb{R}^n, $ $ y\in\Omega.$
	\item Let $\alpha, \beta$ are real, then $$I_{\alpha,\beta}(t) = \int_{\Omega}(\Delta^{\alpha}(y+t))(\Delta^{\beta}(y)) dy<\infty,$$
	if $\beta>-1, \alpha+\beta< 1-\frac{2n}{r}$, and $$I_{\alpha, \beta}(t) = (c_{\alpha, \beta})\Delta^{\alpha+\beta+\frac{n}{r}}(t).$$
	Moreover $$I_{\alpha}(y) = \int_{\mathbb{R}^n}\left|\Delta^{-\alpha}\left(\frac{x+iy}{i}\right)\right| dx<\infty,$$ if $\alpha>\frac{2n}{r}-1,$ and
	$$I_{\alpha}(y) = (c_{\alpha})\Delta^{\alpha+\frac{n}{r}}(y),$$ where $y\in\Omega.$
\end{enumerate}
}
\textbf{Lemma 6.} {\it For any analytic function from $A_{\alpha}^2(T_{\Omega})$ the following integral formula is
valid $$f(z) = \tilde{c}_{\alpha}\int_{T_{\Omega}}B_{\alpha}(z,w)f(w)\, d\nu_{\alpha}(w), z\in T_{\Omega}.\eqno(1)$$}

Let $1\le p<\infty, 1\le q<\infty$, $\frac{n}{r}\le p_1, \frac{1}{p_1}+\frac{1}{p}=1, \frac{n}{r}<\gamma$. Let $f\in A_{\gamma}^{p,q}$, then (1) with $\alpha > \frac{n}{r} - 1$ is valid (Bergman representation formula with $\alpha$ index is valid).

We now collect a few facts on the (possibly weighted) $L^p$-norms of the Bergman kernel and the normalized Bergman kernel. The first result is classical (see, for example, \cite{10}, \cite{14}, \cite{7}).

\textbf{Proposition 1} (Forelly-Rudin estimates). {\it Let $T_{\Omega}\subset\mathbb{C}^n$ be a tubular domains over symmetric cones, and let
$z_0\in T_{\Omega}$ and $1\le p < \infty$. Then $$\int\limits_{T_{\Omega}}|B(\zeta,z_0)|^p\delta^{\beta}(\zeta)d\nu(\zeta) \le
C\delta^{\beta-2(2n/r)(p-1)}(z_0), -1<\beta<(2n/r)(p-1).$$}

The same result is valid for weighted Bergman kernel (see \cite{4}).

We define new Banach mixed norm analytic Bergman-type spaces in $T_{\Omega}\times\ldots\times T_{\Omega}$ in product of tubular domains over symmetric cones  as follows. Let $m\ge1$, $p_j\in(1;\infty)$ $\nu_j>\frac{n}{r}-1$,
$$A_{\Vec{\nu}}^{\Vec{p}} = \left\{f\in H(T_{\Omega}^m) = H(T_{\Omega}\times\ldots\times T_{\Omega}) =\right. $$
$$=\left(\int\limits_{T_{\Omega}}\ldots\left( \int\limits_{T_{\Omega}}|f(z_1,\ldots,z_m)|^{p_1}\Delta^{\nu_1-\frac{n}{r}}(y_1)\,dx_1dy_1\right)^{\frac{p_2}{p1}}\ldots
\Delta^{\nu_m-\frac{n}{r}}(y_m)\,dx_mdy_m\right)^{\frac{1}{p_m}}<\infty\}.$$

Replacing $A$ by $L$ as usual we get larger space of measurable functions with the same norms.

Mixed norm analytic function spaces in tube we just defined above may be defined easily when one $p_j$ index is equal to infinity. As in simpler domains and probably our projection theorems may be extended also to such type analytic function spaces and we pose  this as a problem for readers.

Note first for case of polydisk (when $T_{\Omega}$ is a unit disk) or even $T_{\Omega}$ is a unit ball in $\mathbb{C}^n$ these analytic spaces are not new. They were introduced and studied in \cite{1}-\cite{2}. Note also very similar spaces in $\mathbb{R}^n$ were introduced and studied before by various authors. Our theorem for mentioned particular cases are not new. They can be seen in \cite{1}: For $m=1$ case our theorem is also known (see \cite{12}-\cite{14}).

\textbf{Lemma 7} (Reproducing formulas, see \cite{14}). {\it For any analytic function from $A_{\alpha}^{2}\left(T_{\Omega}\right)$ the following integral formula is valid
$$
f(z)=\widetilde{c}_{\alpha} \int_{T_{\Omega}} B_{\alpha}(z, w)f(w) d v_{\alpha}(w), z \in T_{\Omega}.\eqno(2)
$$}

Let $1 \leq p<\infty, 1 \leq q<\infty, (n / r) \leq p_{1},$ $ \frac{1}{p_{1}}+\frac{1}{p}=1, $ $\left(\frac{n}{r}-1\right)<\gamma$. Let $f \in A_{\gamma}^{p, q},$ then (2) with $\alpha>\left(\frac{n}{r}-1\right)$ is valid (Bergman representation formula with $\alpha$ index is valid).

We provide below a well-known and important application of $r$-lattices of tubular domains over symmentic cones.

\textbf{Lemma 8} (Atomic decomposition of $A_{\gamma}^{p}$, see \cite{14}). {\it Let $p\ge 1$ and $\nu>\frac{n}{r}-1$. Let $\left\{z_{j}\right\}$ be a $\delta$-lattice in $T_{\Omega}, $ $ \delta \in(0,1), $ $ z_{j}=x_{j}+i y_{j},$ $ z_{j} \in T_{\Omega}, $ $j=1, \ldots, r$. Then
$$
\|f\|_{A_{\nu}^{p}} \asymp \sum_{j}\left|f\left(z_{i}\right)\right|^{p} \Delta^{\nu+\frac{n}{r}}\left(y_{i}\right).
$$
Assume that Bergman projection $P_{\nu}$ is bounded on $A_{\nu}^{p}$ and let $\left\{z_{j}\right\}$ be a $\delta$-lattice in $T_{\Omega}$. 

If $f \in A_{\nu}^{p},$ then
$$
f(z)=\sum_{j} \lambda_{j} B_{\nu}\left(z, z_{j}\right) \Delta^{\nu+\frac{n}{r}}\left(y_{i}\right), z \in T_{\Omega},  \eqno(3)
$$
$$ \sum_{j=1}^{\infty}\left|\lambda_{j}\right|^{p} \Delta^{\nu+\frac{n}{r}}\left(y_{j}\right) \leq c\|f\|_{A_{\nu}^{p}}^{p}. \eqno(4)
$$

If
$$
\sum_{j=1}^{\infty}\left|\lambda_{j}\right|^{p}\Delta^{\nu+\frac{n}{r}}\left(y_{j}\right)<+\infty,
$$
then the "sum with $B_{\nu}$" (3) converges in $A_{\nu}^{p}$ and the reverse to (4) is true also.} 

We mention now several known results on Bergman type projections. The weighted Bergman projection $P_{\nu}$ is the orthogonal projection from the Hilbert space $L_{\nu}^{2}\left(T_{\Omega}\right)$ onto its closed subspace $A_{\nu}^{2}\left(T_{\Omega}\right)$ and it is given by the integral formula
$$
\left(P_{\nu} f\right)(z)=\int_{T_{\Omega}} B_{\nu}(z, w) f(w) \Delta^{\nu-\frac{n}{r}}(\operatorname{Im} w) dv(w),
$$
$z \in T_{\Omega},$ $ \nu>\frac{n}{r}-1$ (see  \cite{4}, \cite{1}).

The $L_{\nu}^{p, q}$ boundedness of the Bergman projection $P_{\nu}$ is still an open problem and has attracted a lot of attention in recent years. Today it is only known that this projection extends to a bounded operator on $L_{\nu}^{p, q}$ for general symmetric cones for the range $1 \leq p<\infty, $ $q_{\nu, p}^{\prime}<q<q_{\nu, p}, $ $ q_{\nu, p}=\min \left\{p, p^{\prime}\right\} q_{\nu}, $ $q_{\nu}=1+\frac{\nu}{n / r-1}$ and $\frac{1}{p}+\frac{1}{p^{\prime}}=1$ (see  \cite{4}, \cite{1}).

The importance can be seen, for example, from the following fact. If $P_{\nu}$ extends to a bounded operator on $L_{\nu}^{p, q},$ then the topological dual space $\left(A_{\nu}^{p, q}\right)^{*}$ of the Bergman space $A_{\nu}^{p, q}$ identifies with $A_{\nu}^{p^{\prime}}, $ $ q^{\prime}$ under the integral pairing
$$
\langle f, g\rangle_{\nu}=\int_{T_{\Omega}} f(z) \overline{g(z)} \Delta^{\nu-\frac{n}{r}}(\operatorname{Im} z) dv(z),
$$
$f \in A_{\nu}^{p, q} ; g \in A_{\nu}^{p^{\prime}, q^{\prime}}$ (see   \cite{4}, \cite{1}). Let

$$
\left(T_{\alpha, \beta, \gamma} f\right)(z)=\Delta^{\alpha}(\operatorname{Im} z) \int_{T_{\Omega}} B_{\gamma}(z, w) f(w) \Delta^{\beta}(\operatorname{Im} w) d v(w),
$$
$$
\left(T_{\alpha, \beta, \gamma}^{+} f\right)(z)=\Delta^{\alpha}(\operatorname{Im} z) \int_{T_{\Omega}}\left|B_{\gamma}(z, w)\right| f(w) \Delta^{\beta}(\operatorname{Im} w) d v(w),
$$
$z \in T_{\Omega}, $ $f \in L^{1}(T_\Omega)$. The following assertions were proved in \cite{1}.

\textbf{Theorem A.}  {\it There are $\nu_{1}=\nu_{1}(\alpha, n, r, q), $ $ \nu_{2}=\nu_{2}(\alpha, n, r, q)$ so that for $1 \leq p, q< \infty, $ $ \nu \in \mathbb{R}, $ $\gamma=\alpha+\beta+\frac{n}{r},$ $ \alpha+\beta>-1,$ then $T_{\alpha, \beta, \gamma}^{+}$ is a bounded operator on $L_{\nu}^{p, q}\left(T_{\Omega}\right)$ for all $\nu \in(\nu_{1},$ $ \nu_{2})$.}

\textbf{Theorem B.}  {\it  Let ( $Q^{+}$) be ( $T_{\alpha, \beta, \gamma}^{+}$) operator for $\alpha=0,$ $ \gamma=\nu+m, $ $ \beta=\nu-\frac{n}{r}$. Then $\left(Q^{+}\right)$ for $\nu+m>\frac{n}{r}-1, $ $1 \leq p, q<\infty,$ is a bounded operator from $L_{\nu}^{p, q}$ to $L_{\nu+m q}^{p, q},$ if $\nu \in\left(\nu_{1}, \nu_{2}\right)$ for some $\nu_{1}=\nu_{1}(p, q, n, r, \nu), $ $ \nu_{2}=\nu_{2}(p, q, n, r, \nu), $ $\left(T_{\alpha, \beta, \gamma}^{+}\right)$ is a bounded operator on $L^{\infty},$ if $\alpha>\frac{n}{r}-1, $ $ \beta>-1, $ $\gamma=\alpha+\beta+\frac{n}{r}$. The same is valid for $T_{\alpha, \beta, \gamma}$ operator, acting between same $L^{p,q}$  type spaces into certain $A^{p,q}_\alpha$ spaces.}

It is natural to try to extend these results to larger spaces of analytic function in tube. It is one of the goals of this paper.

Theorems A, B can be seen as direct extensions of well-known old classical theorems on Bergman projection in analytic function spaces in simpler domains to general tube domains over symmetric cones.

\section{\textbf {Main results}}

In this section we formulate main results of this note. Complete and not difficult proofs of these interesting assertions will be provided elsewhere.

\textbf{Theorem 1.} {\it Let
$$T_{\Vec{\beta}}f(\Vec{z}) = \int\limits_{T_{\Omega}^m}\frac{f(w_1,\ldots,w_m)\prod\limits_{j=1}^{m}\Delta^{\beta_j-\frac{n}{r}}(w_j)dv(w_j)}{\Delta^{\beta_1+\frac{n}{r}}(\frac{z_1-\overline{w}_1}{i})\ldots\Delta^{\beta_m+\frac{n}{r}}(\frac{z_m-\overline{w}_m}{i})},$$
$dv(w) = dudv; w = u+iv\in T_{\Omega}$, $\Vec{z} = (z_1,\ldots, z_m)\in T_{\Omega}$. Let $\beta_j>\beta_0, j=1,\ldots,m$ for some fixed enough large $\beta_0$.
Then $T_{\Vec{\beta}}$ operator maps $L_{\Vec{\nu}}^{\Vec{p}}(T_{\Omega}^m)$ into $A_{\Vec{\nu}}^{\Vec{p}}(T_{\Omega}^m)$, $p_j > 1$; $\nu_j>\frac{n}{r}-1$, $j=1,\ldots,m$.}

\textbf{Remark 1.} For unit ball and unit disk this theorem can be seen in \cite{1}, \cite{2}. We provide for simplisity our proof in the unit disk case since repetition of same arguments leads to the proof of theorem 1. The proof use only Minkowski and Holder's inequality and Forelly-Rudin estimate (5) which is avialable in tubular domaines over symmetric cones $\tau > -1$, $\tau_1 > \tau + \frac{2n}{r},$
$$\int\limits_{T_{\Omega}}\frac{\Delta^{\tau}(Im\,w)\cdot dv(w)}{\Delta(Im^{\tau_1}(\frac{w-z}{i}))}\le c\Delta^{\tau-\tau_1+\frac{2n}{r}}(Im z), z\in T_{\Omega}. \eqno(5)$$

\textbf{Remark 2.} In the unit disk for $m = 1$ this result is classical and well-known fact (Bergman projection theorem in tubular domain (\cite{6}, \cite{5})). The proof uses only Forelly-Rudin estimate from Lemma 1 and Holders and Minkowski inequalities and $m = 2$ and unit disk case is typical. We have in the unit disk $U = \{ | z| < 1\}, $ $m = 2$ case the following estimates.

The following theorem for $m = 1$ is well-knowns.

For any two $n$-tuples of real numbers $x = (x_1,\ldots, x_n)$ and $y = (y_1,\ldots, y_n)$ we consider integral operator
$$(R_{x,y}g)(w) = \Delta(Im\,w)^{-m\left(\frac{2n}{r}\right)+\sum\limits_{i=1}^{m}y_i}\,\times
\int\limits_{T_{\Omega}}\ldots\int\limits_{T_{\Omega}}g(z_1,\ldots,z_m)\cdot$$
$$\cdot\prod\limits_{j=1}^{m}\frac{(\Delta(Im\,z_j))^{x_j}}{\Delta(Im(\frac{\overline{z}_j-w}{i}))^{x_j+y_j}}\times dV(z_1)\ldots dV(z_m)$$
for $g\in L^1(T_{\Omega}^m; $ $ dV_{x_1},\ldots, $ $dV_{x_m}); $ $w\in T_{\Omega}, $ $x_j>-1; $ $x_j+y_j > 0; $ $j=1,\ldots,m$.

\textbf{Theorem 2.} {\it Let $s_j > (-1)$ and $s_j < p,$ $ms_j+1 > m(\frac{2n}{r}-y_j) - (m-1)(\frac{2n}{r}),$ $j=1,\ldots,m.$ Then there is exist a constant $C>0$ such that
$$\int\limits_{T_{\Omega}}|R_{x,y}g)(w)|\cdot\Delta(Im\,w)^{(m-1)\frac{2n}{r}+\sum\limits_{j=1}^{m}s_j}\cdot dV(w)\le$$
$$\le C\int\limits_{T_{\Omega}}\ldots\int\limits_{T_{\Omega}}g(z_1,\ldots,z_m)\cdot \prod\limits_{j=1}^{m}\left(\Delta^{s_j}(Im\,z_j)\right)\,dV(z_j).$$}

\textbf{Remark 3.} This theorem is valid probably for all $p>1$. Our theorem 2 for $m=1$ case can be seen in \cite{10} and \cite{12}.

We refer to \cite{12}--\cite{14} for definitions of weighted Hardy $H^1_{\beta}$ and $A^{\infty}_{\beta}$, $\beta\ge 0$ (analytic  Bloch and Hardy type spaces in tube). We formulate a new sharp decomposition theorem for Bergman spaces in tube.

Let $D$ be bounded (or unbounded) domain in $\mathbb{C}^n.$
We seek equivalent relations of the following type
$\|f_1..f_m\|_X \asymp \prod_{j=1}^m \|f_j\|_{X_j},$ where $X,$ $X_j$ are certain classes of analytic function spaces in $D$ domain and $f_1,...,f_m$ are concrete analytic functions from these domains. Such type sharp results were proved by author in various domains and various complex function spaces previously in \cite{20}-\cite{23}. 
In the following theorem we add a new result in this direction. A valuable remark will be added after theorem also.

\textbf{Theorem 3.} {\it We have
$$
\|f_1...f_m\|_{A^1_\alpha}\asymp \prod\|f_j\|_{A^1_{\alpha_k}},$$
for some $\alpha,$ $\alpha_k,$ $k=1,...,m,$ $\alpha >-1,$ $\alpha_k>-1,$ $k=1,...,m,$\\
or
$$
\|f_1...f_m\|_{A^1_{\alpha}(T_{\Omega})}\asymp\prod_{j=1}^{m_0} \|f_j\|_{A^1_{\alpha_j}(T_\Omega)}\times \prod_{j=m_0+1}^{m} \|f_j\|_{X_j},
$$
if
$$
\prod_{i=1}^{m}f_i(\omega_i)=c_{\beta}\int_{T_{\Omega}} \prod_{j=1}^{m}f_j(z)\frac{\Delta^{\beta}(z)d\nu(z)}{\prod_{j=1}^{m}\Delta^{\frac{\beta+2n/r}{m}}(\frac{z-\omega_j}{i})},
$$
where $X_j = A^\infty_{\beta_j}$ or $H^1_{\beta_j},$ $\beta>\beta_0,$ $\beta_0$ is large enough,
$c_\beta$ is a constant of Bergman representation formula, and $w_j \in T_\Omega$ for all $j=1,...,m,$ $m_0\ge 1,$ some fixed indexes $\beta_k\ge 0,$ $k=m_0+1,...,n.$}

Note that if the amount of functions equal to one, then integral representation vanishes and we get a trivial relation.

This type sharp decomposition  theorems in the ball and bounded strongly pseudoconvex domains with smooth boundary were proved earlier in \cite{20}-\cite{23}.

\textbf{Remark 4.} The core of the proof of theorem 3 is a new special integral representation, uniform estimates for Bergman spaces and Forelly-Rudin formula in tube domains over symmetric cones. Probably this sharp decomposition theorem 3 can be extended from $A^1_\alpha$ Bergman analytic function spaces in tubular domains over symmetric cones to all $A^p_\alpha$ spaces in tubular domains for all positive values of $p$ and all $\alpha>-1.$
We leave this interesting question as a problem for interested readers. Note that for strongly pseudoconvex bounded  domains this is true and it was proved previously by author in a recent paper with E. Tomashevskaya (see, for example, \cite{20}-\cite{22}). We remark for interested readers in addition that this sharp decomposition theorem 3 is valid with the same proof for $A^1_\alpha$ Bergman function spaces in very general Siegel domains of second type, and also for bounded symmetric domains.
Another interesting problem is to extend this sharp decomposition theorem 3 from $A^1_\alpha$ function spaces in tubular domains  to general  $A^{p,q}_\alpha$ function spaces for all positive values $p$ and $q$ and and $\alpha>-1$, or at least with some restrictions on positive parameters $p$ and $q$. We leave this rather interesting question  also to interested readers. This type decomposition problems may be considered and solved by similar methods also for Bergman type $A^p_\alpha$
harmonic function spaces of several variables in the upper half space and in the unit ball of $R^n.$ We leave this also to interested readers.

\textbf{Remark 5.} In various other domains this theorem 3 was proved previously by author. We refer the reader to \cite{20}-\cite{23} for similar type results in other analytic function spaces of several variables in other domains in $\mathbb{C}^n.$

All elementar tools (in particular uniform estimates and integral representations and Forelly-Rudin type estimates), which are needed to get similar type new sharp decomposition theorems in Bergman type function spaces of several variables in the unit polydisk and in harmonic function spaces in the unit ball and upper-half spaces can be seen in \cite{30}-\cite{32}. We pose this as a  problem and leave this to various interested readers. We also note that similar type sharp decomposition theorems with very similar proofs may be valid also in minimal bounded homogeneous domains, Siegel domains of second type and also in various bounded symmetric domains in $\mathbb{C}^n.$
 
We refer to \cite{30}-\cite{32} for basic definitions of function theory in the polydisk and unit ball of $\mathbb{R}^n$ of analytic and harmonic functions of several variables. We denote as usual $A^p_\alpha$ Bergman spaces in the polydisk $U^n$ in $\mathbb{C}^n$ and unit ball $B^n$ in $\mathbb{R}^n.$ 

We have similarly as in tube domains in theorem 3 under certain integral condition in polydisk $U^n$ $\|f_1...f_m\|_{A^1_\alpha(U^n)}$ is equivalent to $\prod_{k=1}^m \|f_k\|_{A^1_{\alpha_k}(U^n)}$ 
(under certain Bergman type condition which vanishes for $m=1$) for some $\alpha,$ $\alpha_k,$ $k=1,...,m,$ $\alpha>-1,$ $\alpha_k>-1,$ $k=1,...,m.$ Note that
some modifications of relation which was provided above are also valid if we replace $A^1_{\alpha_k}$ by $A^\infty_{\alpha_k}$ or $H^1_{\alpha_k}.$ Where these are analytic Bloch type and weighted Hardy classes in the polydisk $U^n$ (see \cite{30}-\cite{32} for these spaces in the polydisk), for some $\alpha_k>-1$ (see for such results also theorem 3). 

Very similarly for $A^p_\alpha(B^n)$ harmonic Bergman spaces in the ball $B^n$ in $\mathbb{R}^n$ we have under certain natural integral condition (which vanishes for $m=1$) that $\|f_1...f_m\|_{A^1_\alpha(B^n)}$ is equivalent to
$\prod_{k=1}^m \|f_k\|_{A^1_{\alpha_k}},$ for some $\alpha>-1,$ $\alpha_k>-1,$ $k=1,...,m.$ With obvious simple modifications of the right side of this sharp relation for weighted Hardy and Bloch spaces of harmonic functions spaces of several variables in the ball $B^n$ of $\mathbb{R}^n$.

These assertions can be also shown for $A^1_\alpha$ Bergman harmonic function spaces in $\mathbb{R}^{n+1}$ and related Bloch and weighed Hardy spaces in these unbounded domains
(see \cite{30}-\cite{32} and various references there for these spaces of harmonic functions of several variables in these domains).

Proofs in all cases are very similar and not difficult and can be probably extended even to all values of positive $p,$ namely for $A^p_\alpha$ Bergman spaces in these domains, where $\alpha>-1.$

\section{\textbf {Conclusion}}
The theory of analytic function spaces in rather complicated and general tubular domains over symmetric cones. Is a new research area and our new results on Bergman type projections in these analytic  function spaces may have various interesting  applications in this new research area. In much simpler domains for example such as the unit ball, the unit polydisk, and upper half plane such type interesting  applications are well known in literature. 

All our results have rather transparent proofs and they with the same proof are also valid in various other domains (and various similar type analytic function  spaces on them) such as, for example, Siegel domains of second type or bounded symmetric domains or minimal bounded homogeneous domains and they may have also many applications in complex function theory in these type domains. We leave this to interested readers. These all our results have with the same proof also complete analogues in Bergman type harmonic function spaces of several variables in upperhalf plane and in the unit ball in $\mathbb{R}^n.$ We leave this task also to interested readers.

All our results can be probably extended to mixed norm $A^{pq}_{\alpha}$ spaces in tube and product of tube domains also and to analytic Herz spaces also in tube and in products of tube domains. We refer to \cite{14} for many properties of $A^{pq}_{\alpha}$ mixed norm analytic function spaces in tube. To define such spaces in product domains is an easy task. We define analytic Herz spaces in tube as follows they have the following finite quazinorm
$$
\int_{T_\Omega}\left(\int_{B(w,r)} |f(z)|^p dv_\alpha(z)\right)^{q/p}dv(w),
$$
where $p,q\ge 1,$ $\alpha>-1.$ Such type analytic function spaces in other domains are well studied. 

Or we can replace integration by $T_\Omega$ by summation, and the inner integral by $D(a_k,r),$ where $dv_\alpha$ is weighted Lebegues measure (see \cite{14}) and $(a_k)$ is $r$-lattice in tube (see \cite{14}), and where $B(z,r)$ is a Bergman ball in tube domains, $z\in T_\Omega$ (see \cite{14}).

To define such analytic mixed norm Herz spaces in product domains is an easy task. We pose an extension of all our results of this note to this type function spaces in tube and product of tube domains as a problem to interested readers.

In \cite{33} new Bergman projection theorems were provided for new mixed norm spaces in the polydisk extending well known classical results. It will be interesting to find their complete analogues in tubular domain. We can easily define such spaces also in tube by simple iteration of a norm of $A^{pq}_{\alpha}$ spaces in $T_\Omega$ which we defined above. 

Namely can we say that there is a bounded projection of Bergman type $T_\beta,$ where $\beta$ is large enough, from
$L^{p_1,...,p_m}_{\alpha_1,...,\alpha_m}$ to
$A^{p_1,...,p_m}_{\alpha_1,...,\alpha_m}$ spaces in $T^m_{\Omega}$ products of tube domains, $m\ge 1.$

Where analytic function spaces we indicated above have the following finite quazinorms $\| ...\|f\|_{A^{p_1,p_2}_{\alpha_1}}...\|_{A^{p_{m-1},p_m}_{\alpha_m}},$ $p_j\ge 1,$ $j =1,...,m, $ $\alpha_j>-1,$ $j=1,...,m$ (or related classes of measurable functions with $L$ instead of $A$).

We denote new analytic Herz type spaces defined above by $H^{p,q}_{\alpha},$ $p,q\ge 1,$ $\alpha>-1,$ larger classes consisting from measurable functions by $L^{p,q}_{\alpha},$ for same parameters.
The interesting question is $T_\beta,$ for large enough $\beta,$ that is Bergman projection bounded from spaces with norm $\|...\|f\|_{L^{p_1,p_2}_{\alpha_1}}...\|_{L^{p_{n-1},p_n}_{\alpha_n}}$  to it is analytic subspace
for all $p_j\ge 1,$ $j =1,...,n, $ $\alpha_j>-1,$ $j=1,...,n.$

Let $A^{p,q}_{\alpha}(T_\Omega\times...\times T_\Omega)$
be the space of analytic functions, so that
$$\int_\Omega...\int_\Omega\left(\int_V...\int_V |f|^p dx_1...dx_m\right)^{q/p}\Delta^{\beta_1}(y_1)...
\Delta^{\beta_m}(y_m)dy_1...dy_m$$
is finite. We define large spaces of measurable functions similarly as usual replacing $A$ by $L$. The open interesting question is the following, is there a bounded Bergman projection $T_\beta$ for large enough $\beta$ between these spaces for all $p,q\ge 1,$ and $\beta_j>-1,$ $j=1,...,m.$

We refer the reader to \cite{2}, \cite{40}--\cite{43} for other new interesting results on Bergman projection in various analytic spaces and various domains.

Bergman projection theorems in $\mathbb C^n$ may have many applications in function theory, for example to prove various embedding theorems in $\mathbb C^n$ (see, for example, \cite{15}, \cite{16}, \cite{4}--\cite{3}).

The author thanks Dr. Natalia Makhina for technical help and support.

 \end{document}